\documentclass[twocolumn]{autart}    % Enable this line and disable the
\usepackage{amsmath, amsfonts, amssymb, graphicx, gensymb, textcomp}

\renewcommand{\ll}{\left(\!\!}
\newcommand{\rr}{\!\!\right)}
\newcommand{\Ri}{\mathcal R_i}
\newcommand{\Rb}{\mathcal R_b}
\newcommand{\Rm}{\mathcal R_m}
\newcommand{\Rs}{\mathcal R_s}
\newcommand{\Rmb}{R_{m,b}}
\newcommand{\Rsb}{R_{s,b}}
\renewcommand{\o}{\omega}
\newcommand{\ho}{\hat \o}
\newcommand{\om}{\o_{\max}}
\renewcommand{\to}{\tilde \o}
\newcommand{\cro}[1]{[#1_{\times}]}
\renewcommand{\a}{a}
\renewcommand{\b}{b}
\newcommand{\ao}{\mathring{\a}}
\newcommand{\bo}{\mathring{\b}}
\newcommand{\at}{\tilde \a}
\newcommand{\bt}{\tilde \b}
\newcommand{\aTb}{\a^T\b}
\newcommand{\aTbo}{\ao^T\bo}
\newcommand{\R}[1]{\mathbb R ^{#1}}
\newcommand{\klim}{k^*}
\newcommand{\am}{A_{m}}
\newcommand{\lM}{\lambda_{\max}}
\newcommand{\lm}{\lambda_{\min}}
\newcommand{\ai}[1]{\alpha_{#1}}
\newtheorem{Problem}{Problem}

\newtheorem{assump}{Assumption}
\newtheorem{proposition}{Proposition}

\newtheorem{theorem}{Theorem}
\newtheorem{remark}{Remark}
\begin{document}

\begin{frontmatter}

\title{Angular velocity nonlinear observer from vector measurements}
\thanks[footnoteinfo]{This paper was not presented at any IFAC
meeting. Corresponding author N. Petit. Tel. +33 1 40 51 93 30.
}

\author[MINES]{Lionel Magnis}\ead{lionel.magnis@mines-paristech.fr},    % Add the
\author[MINES]{Nicolas Petit}\ead{nicolas.petit@mines-paristech.fr},               % e-mail address

\address[MINES]{MINES ParisTech, PSL Research University, CAS, 60 bd Saint-Michel, 75272 Paris Cedex FRANCE}  % Please supply

                             % keyword list or with the

%%%%%%%%%%%%%%%%%%%%%%%%%%%%%%%%%%%%%%%%%%%%%%%%%%%%%%%%%%%%%%%%%%%%%%%%%%%%%%%%
\begin{abstract}
The paper proposes a technique to estimate the angular velocity of a rigid body from vector measurements. Compared to the approaches presented in the literature, it does not use attitude information nor rate gyros as inputs. Instead, vector measurements are directly filtered through a nonlinear observer estimating the angular velocity. Convergence is established using a detailed analysis of the linear-time varying dynamics appearing in the estimation error equation. This equation stems from the classic Euler equations and measurement equations. A high gain design allows to establish local uniform exponential convergence. Simulation results are provided to illustrate the method.
\end{abstract}
\end{frontmatter}

%%%%%%%%%%%%%%%%%%%%%%%%%%%%%%%%%%%%%%%%%%%%%%%%%%%%%%%%%%%%%%%%%%%%%%%%%%%%%%%%
\section{Introduction}
This article considers the question of estimating the angular velocity of a rigid body from embedded sensors. This general question is of particular importance in various fields, and in particular for the problem of orientation control. As is well described in~\cite{salcudean1991}, most existing control methods for such second order dynamics require angular velocity information~\cite{boskovic2000,silani2003,lovera2005}. The list of typical control methods employing this information is vast, ranging from Lyapunov control design, feedback linearization, to the computed torque method. Numerous implementations can be found in spacecraft, low-cost unmanned aerial vehicles, guided ammunitions, to name a few.

In the literature, several types of methods have been proposed to address this question. On the one hand, the straightforward solution is to use a strap-down rate gyro~\cite{Titterton2004}, which directly provides measurements of the angular velocities. However, rate gyros being relatively fragile and expensive components, prone to drift, another type of solutions is often preferred. Instead, a \emph{two-step approach} is commonly employed. The \emph{first step} is to determine attitude from vector measurements, i.e. onboard vector measurements of  reference vectors being known in a fixed frame. Vector measurements play a central role in the problem of attitude determination as discussed in a recent survey~\cite{crassidis2007}. In a nutshell, when two independent vectors are measured with vector sensors attached to a rigid body, its attitude can be simply defined as the solution of the classic Wahba problem~\cite{wahba1965} which formulates a minimization problem having the rotation matrix from a fixed frame to the body frame as unknown. The \emph{second step} is to reconstruct angular velocities from the attitude.
At any instant, full attitude information can be obtained~\cite{shuster1978,shuster1990,baritzhack1996,choukroun2003}. In principles, once the attitude is known, angular velocity can be estimated from a time-differentiation. However, noises disturb this process. To address this issue, introducing \emph{a priori} information in the estimation process is a valuable technique to filter-out noise from the estimates. For this reason, numerous observers using the Euler equations for a rigid body have been proposed to estimate angular velocity (or angular momentum, which is equivalent) from full attitude information \cite{salcudean1991,thienel2007,sunde2005,jorgensen2011}.  Besides this two-step approach, a more direct solution can be proposed. In this paper, we expose an algorithm that directly uses the vector measurements and reconstructs the angular velocity in a simple manner.

The contribution of this paper is a nonlinear observer reconstructing the angular velocity of a rotating rigid body from vector measurements \emph{directly}, namely by bypassing the relatively heavy first step of attitude estimation.

The paper is organized as follows. In Section~\ref{sec : notations}, we introduce the notations and the problem statement. We analyze the attitude dynamics (rotation and Euler equations) and relate it to the measurements. In Section~\ref{sec : observer}, we define a nonlinear observer with extended state and output injection. To prove its convergence, the error equation is identified as a linear time-varying (LTV) system perturbed by a linear-quadratic term. The dominant part of the LTV dynamics can be shown, by a scaling resulting from a high gain design, to generate an arbitrarily fast exponentially convergent dynamics. In turn, this property reveals instrumental to conclude on the exponential uniform convergence of the error dynamics. Illustrative simulation results are given in Section~\ref{sec : simu}. Conclusions and perspectives are given in Section~\ref{sec : conclusion}.

\section{Notations and problem statement}
\label{sec : notations}
\subsection{Notations}
\textbf{Vectors in $\R{3}$} are written with small letters $x$. $|x|$ is the Euclidean norm of $x$. $\cro{x}$ is the skew-symmetric cross-product matrix associated with $x$, i.e. ${\forall y \in \R{3}, \ \cro{x}y = x \times y}$. Namely,
\begin{displaymath}
\cro{x} \triangleq \left(
\begin{array}{ccc}
 0   & -x_3 &  x_2 \\
 x_3 &  0   & -x_1 \\
-x_2 &  x_1 &  0
\end{array}
\right)
\end{displaymath}
where $x_1,x_2,x_3$ are the coordinates of $x$ in the standard basis of $\R{3}$.

\textbf{Vectors in $\R{9}$} are written with capital letters $X$. $|X|$ is the Euclidean norm of $X$. The induced norm on $9\times9$ matrices is noted $||\cdot||$. Namely,
\begin{displaymath}
||M|| = \max_{|X|=1} |MX|
\end{displaymath}
For convenience, we may write $X$ under the form
\begin{displaymath}
X = \left(X_1^T,X_2^T,X_3^T\right)^T
\end{displaymath}
with $X_1,X_2,X_3 \in \R{3}$. Note that $${|X|^2 = |X_1|^2+|X_2|^2+|X_3|^2}$$

\textbf{Frames} considered in the following are orthonormal bases of $\R{3}$.
\subsection{Problem statement}
Consider a rigid body rotating with respect to an inertial frame $\Ri$. Note $R$ the rotation matrix from $\Ri$ to a body frame $\Rb$ attached to the rigid body and $\o$ the corresponding angular velocity vector, expressed in $\Rb$.
Assuming that the body rotates under the influence of an external torque $\tau$ (which, is null in the case of free-rotation), the variables $R$ and $\o$ are governed by the following differential equations
\begin{align}
\dot R & = R \cro{\o} \label{eq R}\\
\dot \o & = J^{-1}\left(J\o \times \o + \tau \right) \triangleq E(\o) + J^{-1}\tau \label{eq euler}
\end{align}
where $J = \textrm{diag}(J_1,J_2,J_3)$ is the inertia matrix\footnote{Without restriction, we consider that the axes of $\Rb$ are aligned with the principal axes of inertia of the rigid body.}. Equation~\eqref{eq euler} is known as the set of Euler equations for a rotating rigid body \cite{landau1982}. The torque $\tau$ may result from control inputs or disturbances\footnote{In the case of a satellite e.g., the torque could be generated by inertia wheels, magnetorquers, gravity gradient, among other possibilities.}.

We assume that two reference unit vectors $\ao,\bo$ expressed in $\Ri$ are known, and that sensors arranged on the rigid body allow to measure the corresponding unit vectors expressed in $\Rb$. Namely, the measurements are
\begin{equation}
\label{def : a,b}
\a(t) \triangleq R(t)^T \ao , \quad \b(t) \triangleq R(t)^T \bo
\end{equation}
For implementation, the sensors could be e.g. accelerometers, magnetometers, or Sun sensors to name a few \cite{magnis2014}. We now formulate some assumptions.
\begin{assump}
\label{hypothese ab constants}
$\ao,\bo$ are constant and linearly independent
\end{assump}

\begin{assump}
\label{hypothese tau connu}
$J$ and $\tau$ are known
\end{assump}

\begin{assump}
\label{hypothese O borné}
$\o$ is bounded : $|\o(t)| \leq \om$ at all times
\end{assump}

Assumption~\ref{hypothese ab constants} implies that $${p \triangleq \aTb = \aTbo}$$ is constant for all times. Without loss of generality, we assume $\aTbo \geq 0$ (if not, one can simply consider $-\ao$ instead of $\ao$). The problem we address in this paper is the following.
\begin{Problem}
\label{Problem}
Under Assumptions~\ref{hypothese ab constants}-\ref{hypothese tau connu}-\ref{hypothese O borné}, find an estimate $\ho$ of $\o$ from the measurements $\a,\b$ defined in \eqref{def : a,b}.
\end{Problem}

\section{Observer definition and analysis of convergence}
\label{sec : observer}
\subsection{Observer definition}
The time derivative of the measurement $\a$ is
\begin{equation}
\label{dot ab}
\dot \a = \dot R^T \ao = - \cro{\o} R^T \ao = \a \times \o
\end{equation}
and the same holds for ${\dot \b = \b \times \o}$. To solve Problem~\ref{Problem}, the main idea of the paper is to consider the reconstruction of the extended 9-dimensional state $X$ by its estimate $\hat X$
\begin{displaymath}
X = \left(
\begin{array}{c}
\a \\
\b \\
\o
\end{array}\right), \quad
\hat X = \left(
\begin{array}{c}
\hat a \\
\hat b \\
\ho
\end{array}\right)
\end{displaymath}
The state is governed by
\begin{equation}
\label{eq : X}
\dot X = \left(
\begin{array}{c}
\a \times \o\\
\b \times \o \\
E(\o) + J^{-1}\tau
\end{array}
\right)
\end{equation}
and the following observer is proposed
\begin{equation}
\label{def : observer}
  \dot {\hat X}  =   \left(
\begin{array}{c}
\a \times \ho - \alpha k(\hat \a - \a)\\
\b \times \ho - \alpha k (\hat \b - \b)\\
E(\ho) + J^{-1}\tau +  k^2  \a  \times  (\hat \a - \a) + k^2 \b \times   (\hat \b - \b)
\end{array}
\right)
\end{equation}
where $\alpha \in (0,2\sqrt{1-p})$ and $k>0$ are constant (tuning) parameters. Note
\begin{equation}
\label{def erreur}
\tilde X \triangleq X - \hat X \triangleq \left(
\begin{array}{c}
\at \\
\bt\\
\to
\end{array}\right)
\end{equation}
the error state. We have
\begin{equation}
\label{eq : dot X tilde}
\dot{\tilde X} = \left(
\begin{array}{ccc}
-\alpha k I & 0 & \cro{\a} \\
0 & -\alpha kI & \cro{\b} \\
k^2 \cro{\a} & k^2 \cro{\b} & 0
\end{array}
\right)
  \tilde X
+
\left(
\begin{array}{c}
0 \\
0 \\
E(\o) - E(\ho)
\end{array}
\right)
\end{equation}

In Section~\ref{sec : convergence} we will exhibit, for each value ${\alpha \in (0,2\sqrt{1-p})}$, a threshold value $\klim$ such that for $k > \klim$, $\tilde X$ converges locally uniformly exponentially to zero.

\subsection{Preliminary change of variables and properties}
The study of the dynamics \eqref{eq : dot X tilde} employs a preliminary change of coordinates. Note
\begin{equation}
\label{def : Z}
Z \triangleq \left(
\begin{array}{c}
\tilde \a \\
\tilde \b \\
\frac{\to}{k}
\end{array}\right)
\end{equation}
yielding
\begin{equation}
\label{systeme Z perturbe}
\dot{Z} = kA(t)Z
+
\ll
\begin{array}{c}
0 \\
0 \\
\frac{E(\o) - E(\ho)}{k}
\end{array}
\rr
\end{equation}
with
\begin{equation}
\label{def : A}
A(t) \triangleq \left(
\begin{array}{ccc}
-\alpha I & 0 & \cro{\a(t)} \\
0 & - \alpha I & \cro{\b(t)} \\
\cro{\a(t)} & \cro{\b(t)} & 0
\end{array}
\right)
\end{equation}
which we will analyze as an ideal linear time-varying (LTV) system
\begin{equation}
\label{systeme Z non perturbe}
\dot Z = k A(t) Z
\end{equation}
perturbed by the input term
\begin{equation}
\label{def : xi}
\xi \triangleq \ll
\begin{array}{c}
0 \\
0 \\
\frac{E(\o) - E(\ho)}{k}
\end{array}
\rr
\end{equation}
The idea is that for sufficiently large values of $k$, the rate of convergence of \eqref{systeme Z non perturbe} will ensure stability of system \eqref{systeme Z perturbe}.
We start by upper-bounding $A(t)$ and the disturbance \eqref{def : xi}.
\begin{proposition}[Bound on the unforced LTV system]
\label{prop : A borné}
$A(t)$ defined in \eqref{def : A} is upper-bounded by $$\am \triangleq \max\left(\sqrt{2+2\alpha^2},\sqrt{3+\alpha^2}\right)$$
\end{proposition}
\begin{pf}
Let $Y \in \R{9}$ such that $|Y| = 1$. One has
\begin{align*}
|A(t)Y|^2 = &  |-\alpha Y_1 + \a \times Y_3|^2 + |-\alpha Y_2 + \b \times Y_3|^2 \\
		    &     + |\a \times Y_1 + \b \times Y_2|^2 \\
	   \leq & (1+\alpha^2) \left( |Y_1|^2 + |\a \times Y_3|^2 + |Y_2|^2 + |\b \times Y_3|^2  \right)\\
            & + 2 \left( |\a \times Y_1|^2 + |\b \times Y_2|^2 \right)\\
	   \leq & \max\left(2+2\alpha^2,3+\alpha^2\right) |Y|^2 = \am^2 |Y|^2
\end{align*}
Hence, $||A|| \leq \am$.
\end{pf}

\begin{proposition}[Bound on the disturbance]
For any $Z$, $\xi$ is bounded by
\begin{equation}
\label{prop : borne xi}
|\xi| \leq \sqrt{2}\om |Z| + k |Z|^2
\end{equation}
\end{proposition}
\begin{pf}
We have $$|\xi| = \frac{1}{k} |E(\o) - E(\ho)|$$
with, due to the quadratic nature of $E(\cdot)$,
\begin{align*}
& E(\o) - E(\ho) =  J^{-1} \left( J\to \times \o + J\o \times \to - J\to \times \to \right) \\
& =
\left(
\begin{array}{c}
\frac{J_2-J_3}{J_1} (\o_2 \to_3 + \to_2 \o_3) \\
\frac{J_3-J_1}{J_2} (\o_3 \to_1 + \to_3 \o_1) \\
\frac{J_1-J_2}{J_3} (\o_1 \to_2 + \to_1 \o_2) \\
\end{array}
\right)
-
\left(
\begin{array}{c}
\frac{J_2-J_3}{J_1} \to_2 \to_3 \\
\frac{J_3-J_1}{J_2} \to_3 \to_1 \\
\frac{J_1-J_2}{J_3} \to_1 \to_2 \\
\end{array}
\right) \\
& \triangleq \delta_1 - \delta_2
\end{align*}
As $J_1,J_2,J_3$ are the main moments of inertia of the rigid body, we have \cite{landau1982} (\S 32,9) $$J_i \leq J_j + J_k$$ for all permutations $i,j,k$ and hence
$$\left|\frac{J_2-J_3}{J_1}\right|, \quad \left|\frac{J_3-J_1}{J_2}\right|, \quad \left|\frac{J_1-J_2}{J_3}\right| \quad \leq 1$$
As a straightforward consequence
\begin{displaymath}
|\delta_2| \leq |\to|^2
\end{displaymath}
Moreover, by Cauchy-Schwarz inequality
\begin{displaymath}
(\o_2 \to_3 + \to_2 \o_3)^2 \leq (\o_2^2+\o_3^2) (\to_2^3 + \to_3^2) \leq (\o_2^2+\o_3^2) |\to|^2
\end{displaymath}
Using similar inequalities for all the coordinates of $\delta_1$ yields
\begin{displaymath}
|\delta_1|^2 \leq 2 |\o|^2 |\to|^2 \leq 2 \om^2 |\to|^2
\end{displaymath}
Hence,
\begin{align*}
|\xi| \leq \frac{|\delta_1| + |\delta_2|}{k}&  \leq \sqrt{2} \om \left|\frac{\to}{k}\right| + k \left|\frac{\to}{k}\right|^2 \\
	  & \leq \sqrt{2} \om |Z| + k |Z|^2	
\end{align*}
\end{pf}

\subsection{Analysis of the LTV dynamics $\dot Z = k A(t) Z$}
\label{sec : Z'=kAZ}
We will now use a result on the exponential stability of LTV systems. The claim of \cite{hill2011} Theorem 2.1, which is instrumental in the proof of the next result, is as follows: consider a LTV system ${\dot Z = M(t) Z}$ such that
\begin{itemize}
\item $M(\cdot)$ is $l-$Lipschitz
\item there exists $K\geq 1, c \geq 0$ such that for any $t$ and any $s\geq 0$, ${||e^{M(t)s} || \leq Ke^{-cs}}$
\end{itemize}
Then, for any $t_0, Z_0$, the solution of $\dot Z = M(t)Z$ with initial condition $Z(t_0) = Z_0$ satisfies, for any $t\geq t_0$,
\begin{displaymath}
|Z(t)| \leq K e^{(\sqrt{K l \ln K}-c)(t-t_0)} |Z_0|
\end{displaymath}
Using this result, we will show that the convergence of~\eqref{systeme Z non perturbe} can be tailored by choosing $k$ to arbitrarily increase the rate of convergence, while keeping the overshoot constant.
\begin{theorem}
\label{thm systeme non perturbe}
Let $\alpha \in (0,2\sqrt{1-p})$ be fixed. There exists a continuous function $\gamma(k)$ satisfying
\begin{displaymath}
\lim_{k\rightarrow +\infty} \gamma(k) = + \infty
\end{displaymath}
such that the solution of \eqref{systeme Z non perturbe} satisfies $$|Z(t)| \leq K e^{-\gamma(k)(t-t_0)}|Z(t_0)|$$
with
\begin{equation}
\label{def : K}
K \triangleq \sqrt{\frac{{1+\frac{\alpha}{2\sqrt{1-p}}}}{1-\frac{\alpha}{2\sqrt{1-p}}}}
\end{equation}
for any initial condition $t_0, Z(t_0)$ and any ${t\geq t_0}$.
\end{theorem}
\begin{pf}
Consider any fixed value of $t$. We start by studying the frozen-time matrix $A(t)$. Note
\begin{align*}
\mu & \triangleq \sqrt{8(1-p^2)}
\end{align*}
Introduce the following (time-varying) matrices
\begin{align*}
P_1 & = \left(
\begin{array}{ccc}
\a & 0 & \frac{b-pa}{\sqrt{2(1-p^2)}} \\
0  & b & \frac{a-pb}{\sqrt{2(1-p^2)}}\\
0 & 0 & 0
\end{array}
\right) \\
P_2 & = \frac{1}{\mu} \left(
\begin{array}{cc}
2(p\a-\b) & 0 \\
2(\a-p\b) & 0 \\
\alpha \ \a \times \b & -\sqrt{8-\alpha^2} \ \a \times \b
\end{array}
\right) \\
P_3 &  = \frac{1}{\mu} \left(
\begin{array}{cc}
2 \a \times \b & 0\\
2 \a \times \b & 0\\
\alpha (\b-\a) & \sqrt{4(1+p)-\alpha^2} (\a-\b)
\end{array}
\right) \\
P_4 & = \frac{1}{\mu}\left(
\begin{array}{cc}
2 \a \times \b & 0\\
- 2 \a \times \b & 0\\
\alpha (\a+\b) & - \sqrt{4(1-p)-\alpha^2} (\a+\b)
\end{array}
\right)
\end{align*}
and $$P = (P_1|P_2|P_3|P_4) \in \R{9 \times 9}$$
We have
\begin{displaymath}
P^{-1} A(t) P = \left(\begin{array}{cccc}
M_1 & 0 & 0 & 0 \\
0 & M_2 & 0 & 0 \\
0 & 0 & M_3 & 0 \\
0 & 0 & 0 & M_4
\end{array}\right)
\end{displaymath}
with
\begin{align*}
M_1 & = -\alpha I, \quad M_i = \frac{1}{2}\left(\begin{array}{cc}
-\alpha & -\sqrt{\ai{i}-\alpha^2} \\
\sqrt{\ai{i}-\alpha^2} &  -\alpha
\end{array}\right)
\end{align*}
for $i=2,3,4$ with
\begin{displaymath}
\ai{2} \triangleq 2\sqrt{2} \ > \ \ai{3} \triangleq 2 \sqrt{1+p} \ \geq  \ \ai{4} \triangleq 2\sqrt{1-p} \ >\ \alpha
\end{displaymath}
For all $s \geq 0$
\begin{displaymath}
||e^{A(t)s}|| \leq ||P|| \ ||P^{-1}|| \ e^{-\frac{\alpha}{2}s}
\end{displaymath}
Moreover
\begin{displaymath}
||P|| \ ||P^{-1}|| = \sqrt{\frac{\lM(P^T P)}{\lm(P^TP)}}
\end{displaymath}
where $\lM, \ \lm$ respectively designate the maximum and minimum eigenvalues. Besides,
\begin{displaymath}
P^TP = \left(\begin{array}{cccc}
I & 0 & 0 & 0 \\
0 & Q_2 & 0 & 0 \\
0 & 0 & Q_3 & 0 \\
0 & 0 & 0 & Q_4
\end{array}\right)
\end{displaymath}
with, for $i=2,3,4$
\begin{displaymath}
Q_i  =\left(\begin{array}{cc}
1+ \frac{\alpha^2}{\ai{i}^2} & \frac{\alpha}{\ai{i}} \sqrt{1-\frac{\alpha}{\ai{i}}} \\
\frac{\alpha}{\ai{i}} \sqrt{1-\frac{\alpha}{\ai{i}}} & 1-\frac{\alpha^2}{\ai{i}^2}
\end{array}\right)
\end{displaymath}
yielding the eigenvalues
\begin{displaymath}
\textrm{eig}(P^T P) = \left\{1,1\pm \frac{\alpha}{2\sqrt{2}}, 1\pm \frac{\alpha}{2\sqrt{1+p}}, 1\pm \frac{\alpha}{2\sqrt{1-p}}, \right\}
\end{displaymath}
Thus, for all $s \geq 0$
\begin{displaymath}
||e^{A(t)s}|| \leq K e^{-\frac{\alpha}{2}s}
\end{displaymath}
with
\begin{displaymath}
K = \sqrt{\frac{\lM(P^T P)}{\lm(P^TP)}} = \sqrt{\frac{{1+\frac{\alpha}{2\sqrt{1-p}}}}{1-\frac{\alpha}{2\sqrt{1-p}}}}
\end{displaymath}
Let $k>0$ be fixed. The scaled matrix $kA(\cdot)$ satisfies
$$||e^{kA(t)s}|| \leq K e^{-\frac{k\alpha}{2} s}, \quad \forall t, \ \forall s\geq 0$$
Moreover, for any $Y \in \R{9}$ and any $t,s \in \mathbb R$, one has
\begin{align*}
& (kA(s)-kA(t))Y = k\int_t^s \dot A(u) du \ Y \\
			 & = k\left(
			 \begin{array}{c}
			 \int_t^s \a(u) \times \o(u) du \ Y_3 \\
			 \int_t^s \b(u) \times \o(u) du \ Y_3 \\
			 \int_t^s \a(u) \times \o(u) du \ Y_1 + \int_t^s \b(u) \times \o(u) du \ Y_2 \\
			 \end{array}
			 \right)
\end{align*}
Hence
\begin{align*}
|(kA(s)-kA(t))Y|^2 & \leq  2 \om^2 k^2 |s-t|^2 |Y|^2
\end{align*}
Thus, $kA(\cdot)$ is $kL-$Lipschitz with
\begin{equation}
\label{def : L}
L \triangleq \sqrt{2}\om
\end{equation}
We now apply \cite{hill2011}, Theorem 2.1. For any $t_0$ and any $Z_0$, the solution of \eqref{systeme Z non perturbe} with initial condition $Z(t_0) = Z_0$ satisfies for all $t\geq t_0$
\begin{align*}
|Z(t)| & \leq K e^{(\sqrt{KkL\ln K}-\frac{k\alpha}{2})(t-t_0)} |Z_0|
\end{align*}
which concludes the proof with
\begin{align}
\label{def : gamma}
\gamma(k) & \triangleq \frac{k\alpha}{2} - \sqrt{KkL\ln K}
\end{align}
\end{pf}

\begin{remark} Additionally, one can note that
\begin{itemize}
\item ${\gamma(k) > 0 \quad \Leftrightarrow \quad k > \frac{4K \ln K L}{\alpha^2}}$ \
in which case Theorem~\ref{thm systeme non perturbe} ensures exponential stability of system \eqref{systeme Z non perturbe}.
\item $\gamma(\cdot)$ is strictly increasing for ${k > \frac{4K \ln K L}{\alpha^2}}$.
\end{itemize}
\end{remark}

\subsection{Convergence of the observer}
\label{sec : convergence}
Define $r$ as
\begin{equation}
\label{def : rayon d'attraction}
r(k) \triangleq \frac{1}{\sqrt{\am} K^3} \left(1-\frac{K^2 \sqrt 2 \om}{\gamma(k)} \right) \left( \frac{\gamma(k)}{k} \right)^{\frac{3}{2}}
\end{equation}
and $\klim$ as
\begin{equation}
\label{def : klim}
\klim = \frac{\left(\sqrt{\ln K} + \sqrt{\ln K + 2\alpha K}\right)^2}{\alpha^2} \sqrt 2 K \om > 0
\end{equation}
The following holds
\begin{proposition}
$r(k) > 0$ if and only if $k > \klim$
\end{proposition}
\begin{pf}
A simple rewriting of $r(k)>0$ yields, successively,
\begin{align*}
& r(k) > 0  \Leftrightarrow \gamma(k) > K^2\sqrt 2 \om  = K^2 L\\
& \Leftrightarrow \frac{\alpha}{2} k - \sqrt{LK \ln K} \sqrt{k} - K^2 L > 0 \\
& \Leftrightarrow \sqrt{k} > \frac{\sqrt{LK \ln K} + \sqrt{LK \ln K + 2\alpha LK^2}}{\alpha} = \sqrt{\klim}
\end{align*}
which concludes the proof.
\end{pf}
We can now state the main result of the paper.
\begin{theorem}[main result]
\label{thm convergence obs}
For any ${\alpha \in (0,2\sqrt{1-p})}$, there exists $\klim$ defined by~\eqref{def : klim}-\eqref{def : K} such that for ${k > \klim}$, the observer~\eqref{def : observer} defines an error dynamics~\eqref{eq : dot X tilde} for which the equilibrium 0 is locally uniformly exponentially stable. The basin of attraction of this equilibrium contains the ellipsoid
\begin{equation}
\label{bassin d'attraction}
\left \{\tilde X(0), \quad |\tilde \a (0)|^2 + |\tilde \b (0)|^2 + \frac{|\to(0)|^2}{k^2} < r(k)^2 \right \}
\end{equation}
where $r(k)$ is defined by~\eqref{def : rayon d'attraction}.
\end{theorem}

\begin{pf}
Let $k>k^*$.  Consider the candidate Lyapunov function
\begin{displaymath}
V(t,Z) \triangleq Z^T \left( \int_t^{+\infty} \phi(\tau,t)^T \phi(\tau,t) d\tau \right)Z
\end{displaymath}
where $\phi$ is the transition matrix of system~\eqref{systeme Z non perturbe}. Let $(t,Z)$ be fixed. From Proposition~\ref{prop : A borné}, $kA(\cdot)$ is bounded by $k\am$. Thus (see for example \cite{khalil2000} Theorem 4.12)
\begin{displaymath}
V(t,Z) \geq \frac{1}{2k\am} |Z|^2 \triangleq c_1 |Z|^2 \triangleq W_1(Z)
\end{displaymath}
Moreover, Theorem~\ref{thm systeme non perturbe} implies that for all $\tau \geq t$
\begin{displaymath}
| \phi(\tau,t)Z| \leq K e^{-\gamma(k) (\tau-t)}|Z|
\end{displaymath}
which gives
\begin{align*}
V(t,Z) & \leq K^2 \int_t^{+\infty} e^{-2\gamma(k)(\tau-t)} d\tau |Z|^2 = \frac{K^2}{2\gamma(k)} |Z|^2 \\
       & \triangleq c_2 |Z|^2 \triangleq W_2(Z)
\end{align*}
By construction, $V$ satisfies
\begin{displaymath}
\frac{\partial V}{\partial t}(t,Z) + \frac{\partial }{\partial Z}V(t,Z) kA(t) Z = -|Z|^2
\end{displaymath}
Hence, the derivative of $V$ along the trajectories of \eqref{systeme Z perturbe} is
\begin{align*}
\frac{d}{dt} V(t,Z) & = -|Z|^2 + \frac{\partial V}{\partial Z}(t,Z) \ \xi
\end{align*}
Using
\begin{align*}
\left|\frac{\partial}{\partial Z} V(t,Z)\right| & = 2 \left| \int_t^{+\infty} \phi(\tau,t)^T \phi(\tau,t) d\tau Z \right| \leq \frac{K^2}{\gamma(k)} |Z|
\end{align*}
together with inequality \eqref{prop : borne xi} yields
\begin{displaymath}
\left|\frac{\partial V}{\partial Z}(t,Z) \ \xi \right| \leq \frac{K^2}{\gamma(k)} \left(\sqrt 2 \o_{\max} |Z|^2 + k |Z|^3\right)
\end{displaymath}
Hence
\begin{align*}
\frac{d}{dt} V(t,Z) & \leq - |Z|^2\left(1 - \frac{K^2 \sqrt 2 \om}{\gamma(k)} - \frac{kK^2}{\gamma(k)}|Z|\right) \\
& \triangleq -W_3(Z)
\end{align*}
As $k > \klim$, we have
\begin{align*}
1-\frac{K^2 \sqrt 2 \om}{\gamma(k)} > 0
\end{align*}
We proceed as in \cite{khalil2000} Theorem~4.9. If the initial condition of \eqref{systeme Z perturbe} satisfies
\begin{align*}
& |Z(0)| < r(k) \\
\Leftrightarrow &  |Z(0)| < \frac{\gamma(k)}{k K^2} \left(1-\frac{K^2 \sqrt 2 \om}{\gamma(k)}\right) \times \sqrt{\frac{c_1}{c_2}}
\end{align*}
then $W_3(Z(0)) >0$ and,  while $W_3(Z(t)) >0$, $Z(\cdot)$ remains bounded by
\begin{align*}
|Z(t)|^2 & \leq \frac{V(t)}{c_1} \leq \frac{V(0)}{c_1} \leq \frac{c_2}{c_1} |Z(0)|^2
\end{align*}
which shows that
\begin{align*}
W_3(Z) \geq \left(1-\frac{K^2 \sqrt{2} \om}{\gamma(k)} - \frac{kK^2}{\gamma(k)}\sqrt{\frac{c_2}{c_1}} |Z(0)| \right) |Z|^2
\end{align*}
From \cite{khalil2000}, Theorem~4.10, \eqref{systeme Z perturbe} is locally uniformly exponentially stable. From \eqref{def : Z}, one directly deduces that the basin of attraction contains the ellipsoid \eqref{bassin d'attraction}.
\end{pf}
\begin{remark}
The limitations imposed on $\tilde \a(0)$ and $\tilde \b(0)$ in \eqref{bassin d'attraction} are not truly restrictive, as the actual values $\a(0),\b(0)$ are assumed known, so the observer may be initialized with ${\at(0) = 0, \ \bt(0) = 0}$. What matters is that the error on the unknown quantity $\o(0)$ can be large in practice. Interestingly, when $k$ goes to infinity $r(k)$ tends to the limit
\begin{displaymath}
\frac{1}{\sqrt{\am} K^3} \left(\frac{\alpha}{2}\right)^{\frac{3}{2}} > 0
\end{displaymath}
and arbitrarily large $\to(0)$ is thus allowed from~\eqref{bassin d'attraction}.
\end{remark}
\begin{remark}
The threshold $\klim$ depends linearly on $\om$, which gives helpful hint in the tuning of observer~\eqref{def : observer}.
\end{remark}
\section{Simulation results}
\label{sec : simu}
In this section we illustrate the dependence of the observer with respect to three parameters
\begin{itemize}
\item $p$ which quantifies the linear independence of $(\ao,\bo)$
\item $\om$ the maximal rotation rate of the rigid body
\item the tuning gain $k$
\end{itemize}

Simulations were run for a model of a CubeSat \cite{cubesat2014}. The rotating rigid body under consideration is a rectangular
parallelepiped of dimensions ${20~\textrm{cm} \times 10~\textrm{cm} \times 10~\textrm{cm}}$ and mass $2$kg assumed to be homogeneously distributed. No torque is applied on this system, which is thus in free-rotation.

\begin{figure}[!ht]
\includegraphics[width = \columnwidth]{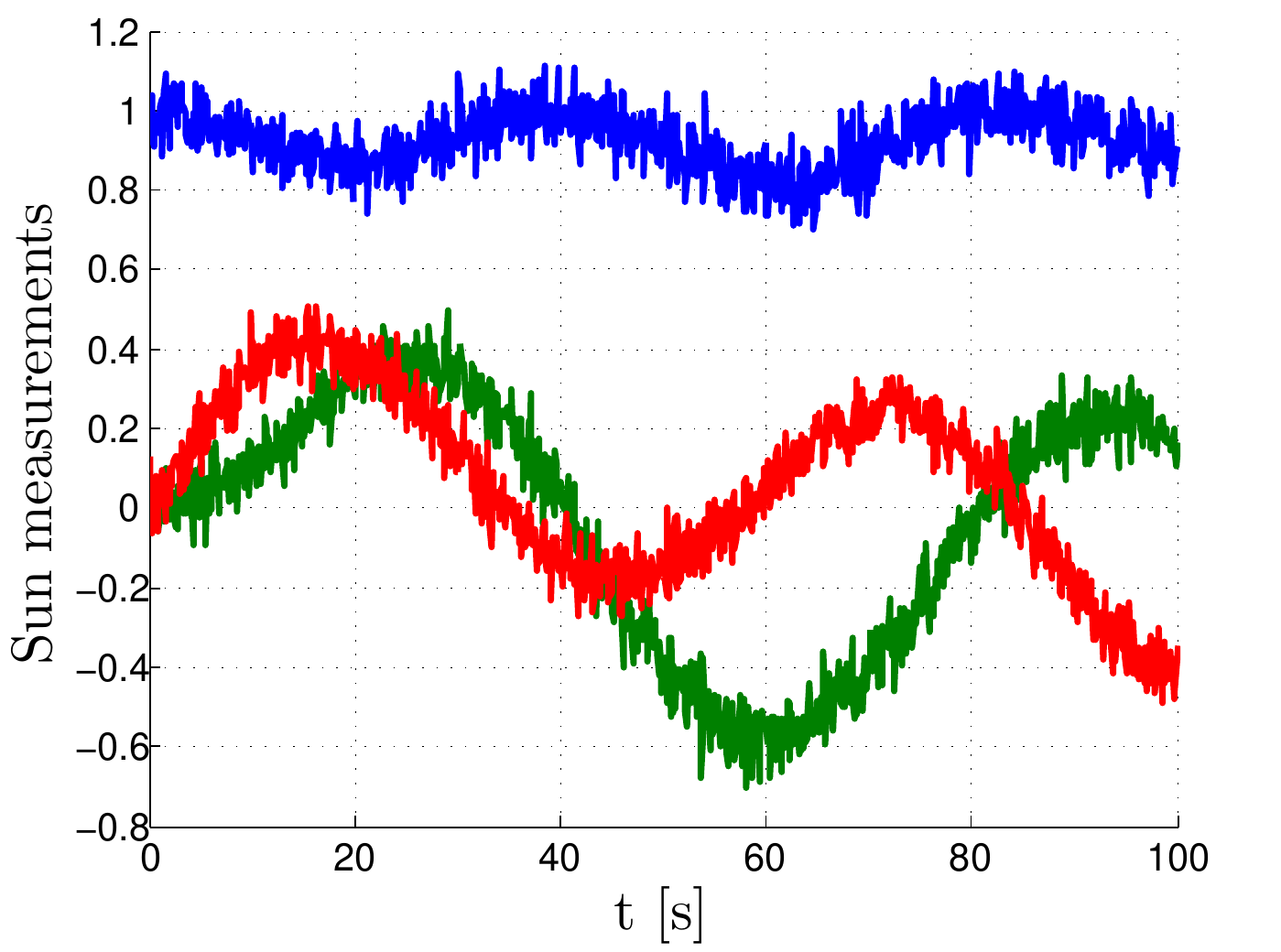}
\includegraphics[width = \columnwidth]{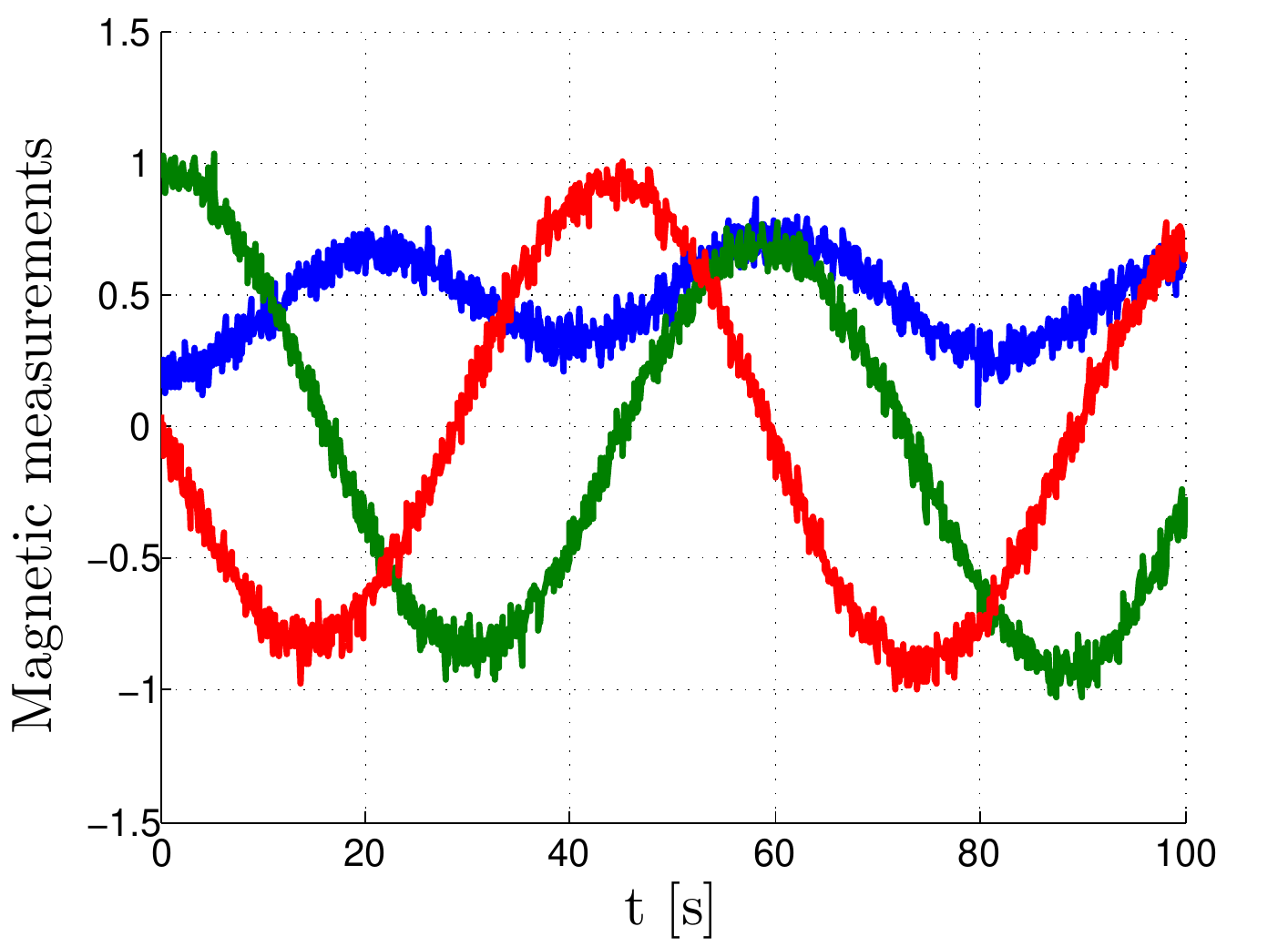}
\caption{Normalized sensor outputs during rotation motion: Sun (top, three coordinates) and magnetic field (bottom, three coordinates)}
\label{fig : mesures}
\end{figure}

In this simulation the two reference unit vectors are the Sun direction $\ao$ and normalized magnetic field $\bo$. The satellite is equipped with
\begin{itemize}
\item 6 Sun sensors providing at all times a measure of the Sun direction $y_a$ in a Sun sensor frame $\Rs$
\item 3 magnetometers able to measure the normalized magnetic field $y_b$ in a magnetometer frame $\Rm$
\end{itemize}
Typical sensor outputs are given in Figure~\ref{fig : mesures}. Because the initial angular velocity vector is not aligned with any of the principal axes of inertia, the rotation motion is not periodic. As can be observed, significant levels of noise have been added on each channel.

It shall be noted that, in practical applications, the sensor frames $\Rs$ need not coincide $\Rm$ and can also differ from the body frame $\Rb$ (defined along the principal axes of inertia) through a constant rotation $\Rmb$, respectively $\Rsb$. With these notations, we have
\begin{displaymath}
a = \Rmb^T y_a, \quad b = \Rsb^T y_b
\end{displaymath}
which is a simple change of coordinates of the measurements.

For sake of accuracy in the implementation, reference dynamics~\eqref{eq : X} and state observer~\eqref{def : observer}
  were simulated using Runge-Kutta 4 method with sample period $0.1$s for various values of $p$ and $\omega(0)$ and with $\alpha = \sqrt{1-p}$.

\begin{figure}[!h]
\includegraphics[width=\columnwidth]{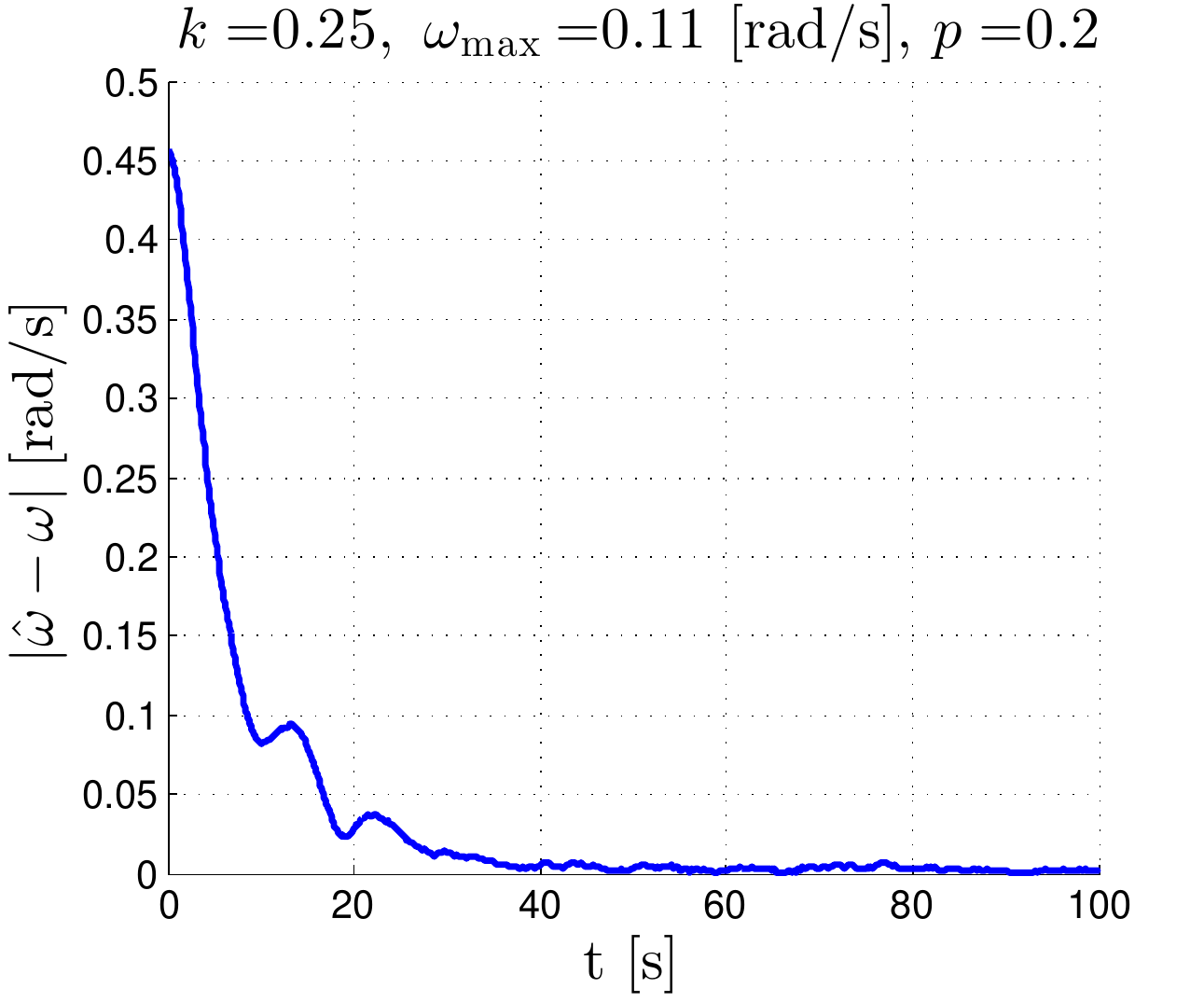}
\caption{Convergence of the observer}
\label{fig : ref}
\end{figure}

\begin{figure}[!h]
\includegraphics[width=\columnwidth]{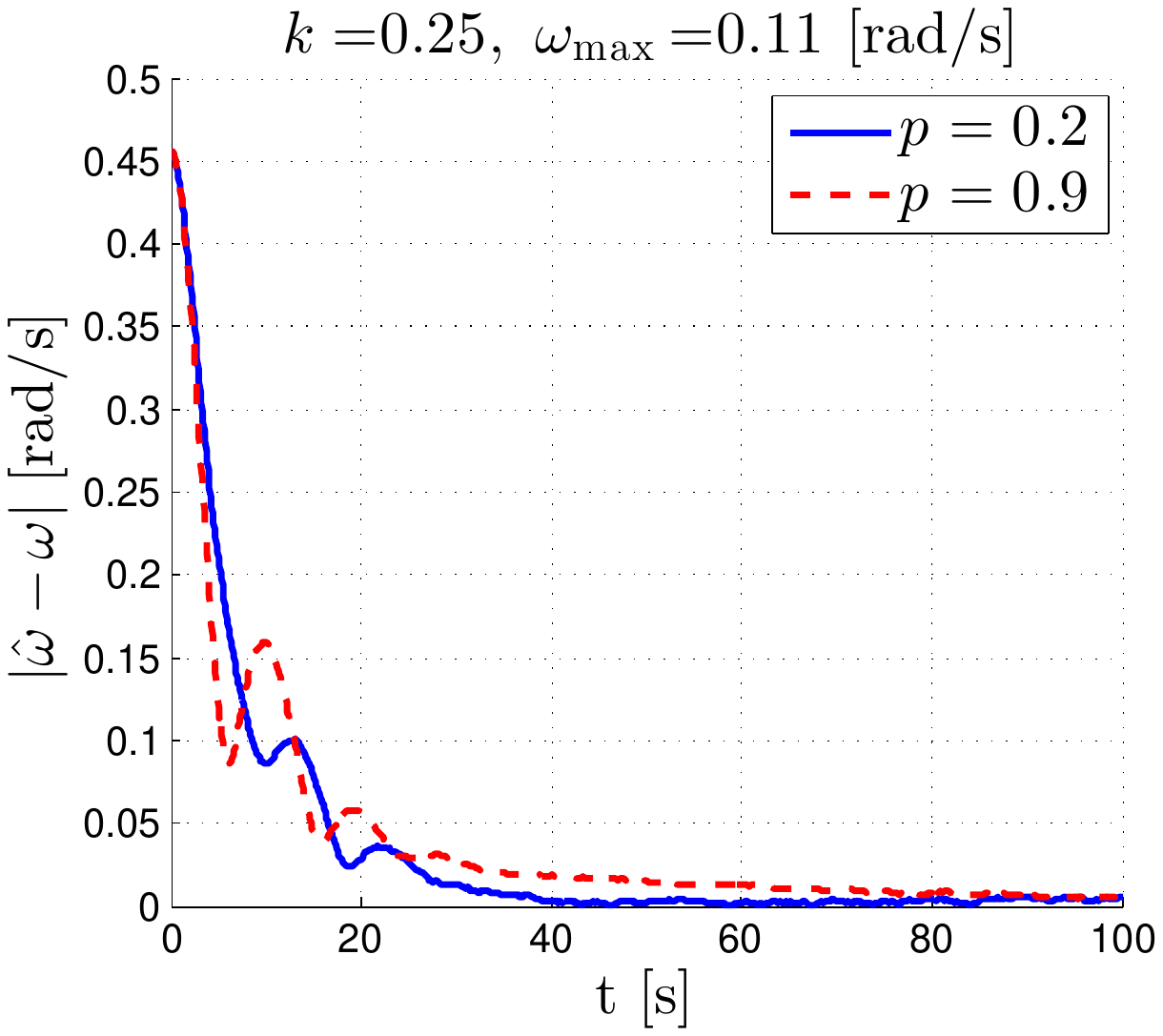}
\caption{The rate of convergence degrades when $p$ increases.}
\label{fig : p}
\end{figure}

\begin{figure}[!h]
\includegraphics[width=\columnwidth]{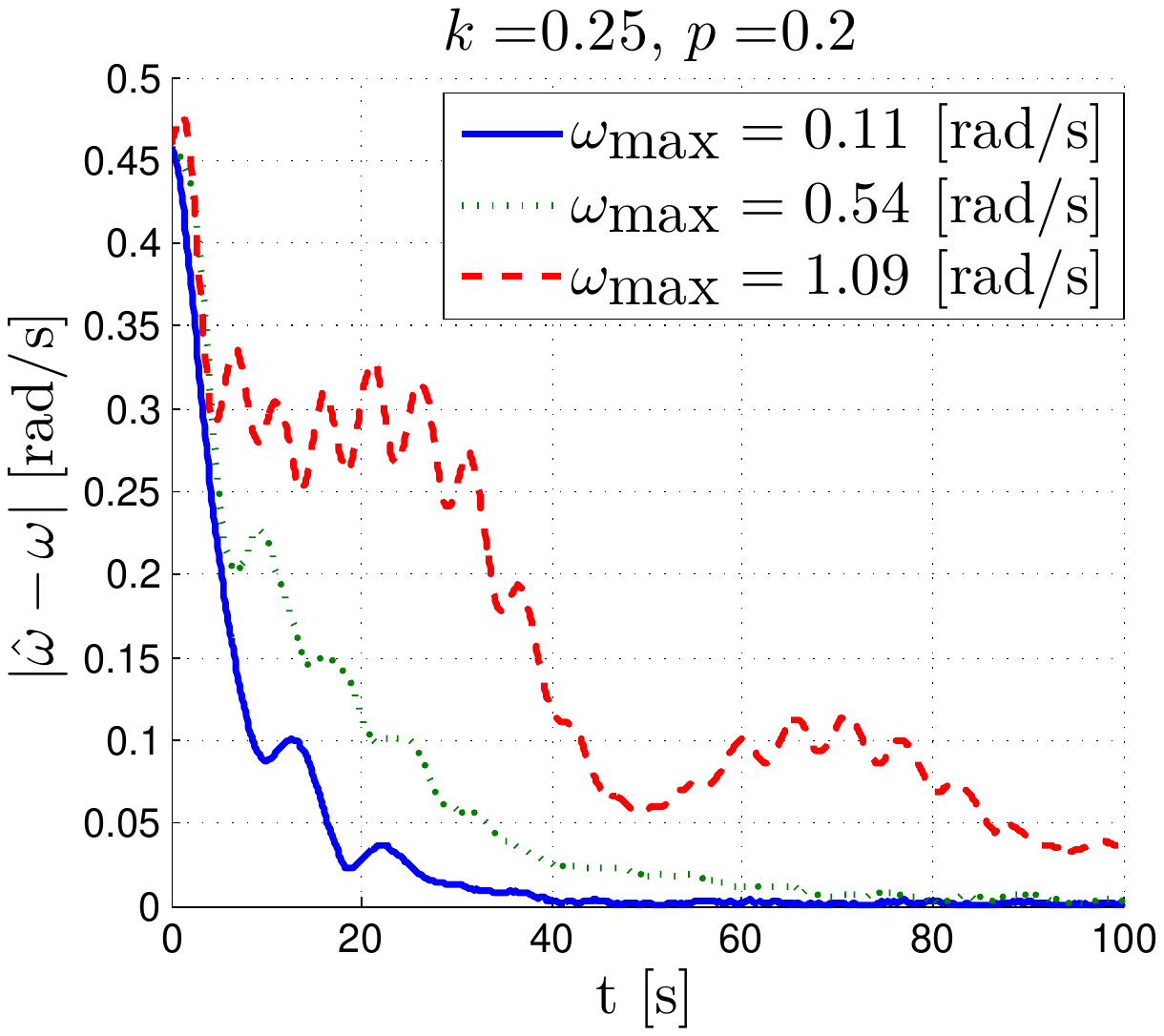}
\caption{Impact of $\om$ on the convergence rate}
\label{fig : omega}
\end{figure}

\begin{figure}[!h]
\includegraphics[width=\columnwidth]{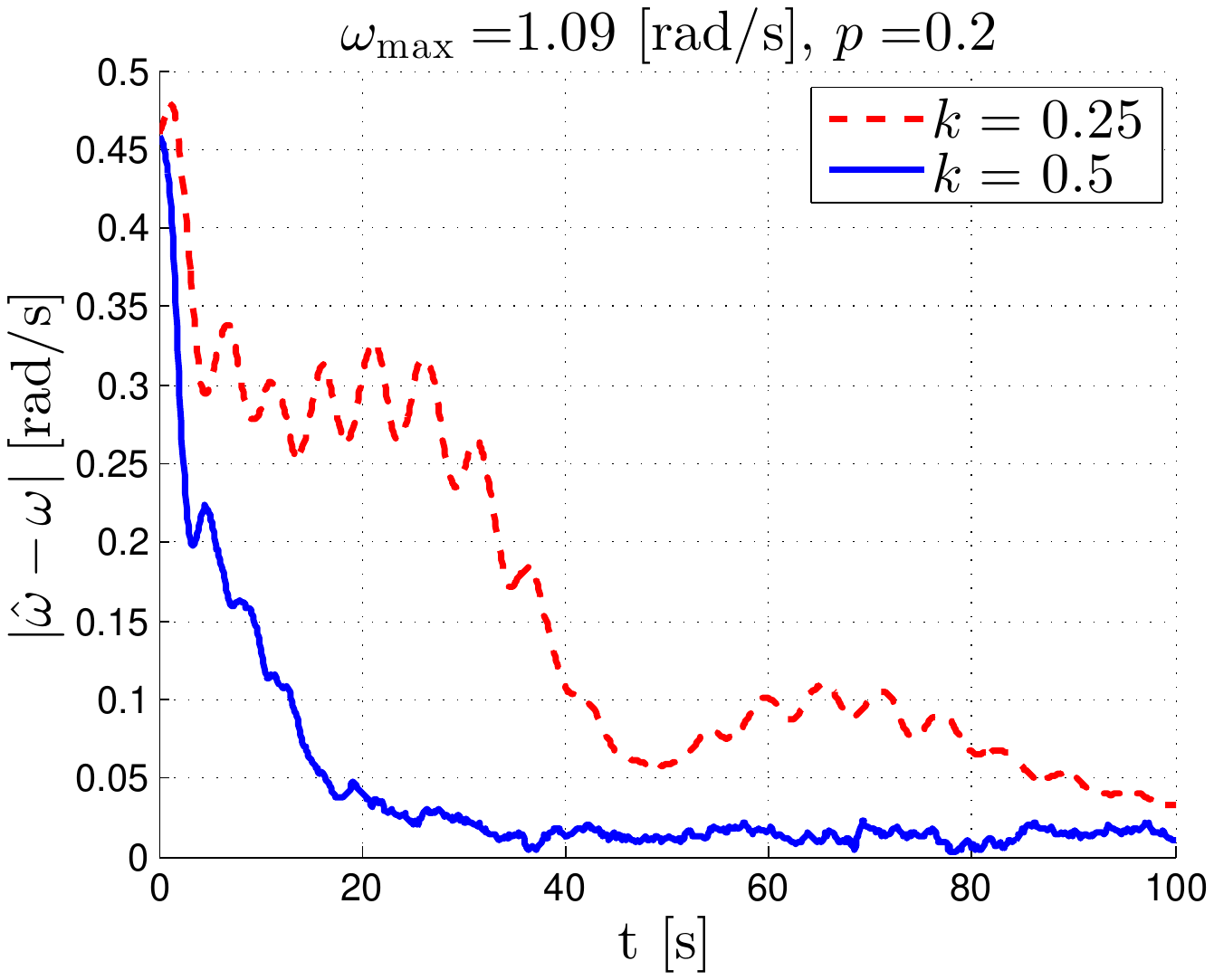}
\caption{When $k$ increases, the convergence is faster but the measurement noise filtering degrades.}
\label{fig : k}
\end{figure}

Figure~\ref{fig : ref} shows the convergence of the observer with parameter values corresponding to the measurements shown in Figure~\ref{fig : mesures}. Note that the vector measurement noise is smoothly filtered by the observer, thanks to the relatively low value of the gain $k$.
Figure~\ref{fig : p} shows the influence of $p$. When $p$ gets close to 1, the rate of convergence is decreased. This was to be expected. To the limit, when $p=1$, all the matrices $A(t)$ become singular and the proof of convergence can not be applied anymore.
In Figure~\ref{fig : omega} we report the behavior of the observer for increasing values of $\om$. The faster the rotation, the slower the convergence. A faster convergence can be achieved by increasing the gain $k$. This increases the sensitivity to noise, as represented in Figure~\ref{fig : k}.
\section{Conclusions and perspectives}
\label{sec : conclusion}
A new method to estimate the angular velocity of a rigid body has been proposed in this article. The method uses onboard measurements of constant and independent vectors. The estimation algorithm is a nonlinear observer which is very simple to implement and induces a very limited computational burden. At this stage, an interesting (but still preliminary) conclusion is that, in the cases considered here, rate gyros could be replaced with an estimation software employing cheap, rugged and resilient sensors. In fact, any set of sensors producing vector measurements such as e.g., Sun sensors, magnetometers, could constitute one such alternative. Assessing the feasibility of this approach requires further investigations including experiments.

More generally, this observer should be considered as a first element of a class of estimation methods which can be developed to address several cases of practical interest. In particular, the introduction of noise in the measurement and uncertainty on the input torque (assumed here to be known) will require extensions such as optimal filtering to treat more general cases. White or colored noises will be good candidates to model these elements. Also, slow variations of the reference vectors $\ao$, $\bo$ should deserve particular care, because such drifts naturally appear in some cases. For example, the Earth magnetic field measured onboard satellites varies according to the position along the orbit.

On the other hand, one can also consider that this method can be useful for other estimation tasks. Among the possibilities are the estimation of the inertia $J$ matrix which we believe is possible from the measurements considered here. This could be of interest for the recently considered task of space debris removal~\cite{bonnal2013}.
Finally, recent attitude estimation techniques have favored the use of vector measurements \emph{together} with rate gyros measurements as inputs.
Among these approaches, one can find {\emph{i)} Extended Kalman Filters (EKF)-like algorithms e.g. \cite{choukroun2006,schmidt2008}, {\emph{ii)} nonlinear observers~\cite{mahony2008,martin2010,vasconcelos2008,tayebi2011,grip2011,trumpf2012}.
This contribution suggests that, here also, the rate gyros could be replaced with more in-depth analysis of the vector measurements.

%\bibliographystyle{unsrt}
%\bibliography{../bibfile}

\end{document}